
\magnification=\magstep0




\font\bigbf=cmbx12 
\font\smc=cmcsc10 
\font\tenjap=msbm10 
\font\sevenjap=msbm7

\def\bbC{{\mathchoice{\hbox{\tenjap C}}{\hbox{\tenjap C}}
{\hbox{\sevenjap C}}{\hbox{\sevenjap C}}}}
\def\bbP{{\mathchoice{\hbox{\tenjap P}}{\hbox{\tenjap P}}
{\hbox{\sevenjap P}}{\hbox{\sevenjap P}}}}
\def\bbR{{\mathchoice{\hbox{\tenjap R}}{\hbox{\tenjap R}}
{\hbox{\sevenjap R}}{\hbox{\sevenjap R}}}}
\def\bbZ{{\mathchoice{\hbox{\tenjap Z}}{\hbox{\tenjap Z}}
{\hbox{\sevenjap Z}}{\hbox{\sevenjap Z}}}}

\def\Kahler{{K\"ahler}}
\def\w{{\mathchoice{\,{\scriptstyle\wedge}\,}{{\scriptstyle\wedge}}
{{\scriptscriptstyle\wedge}}{{\scriptscriptstyle\wedge}}}}


\def\cA{{\cal A}}\def\cH{{\cal H}}\def\cL{{\cal L}}
\def\cO{{\cal O}}\def\cP{{\cal P}}\def\cT{{\cal T}}\def\cU{{\cal U}}


\centerline{\bigbf Some Examples of }
\centerline{\bigbf Special Lagrangian Tori}
\vskip 20 pt
\centerline{by}
\centerline{\smc Robert L. Bryant}
\vskip 10 pt
\centerline{Department of Mathematics}
\centerline{Duke University}
\centerline{Durham, NC 27708}
\vskip 10 pt
\centerline{\tt bryant@math.duke.edu}
\vskip 10 pt
\centerline{12 February 1998}

\bigskip\bigskip
\centerline{\bf \S0. Introduction}
\nobreak\bigskip\nobreak
The purpose of this note is to point out some very elementary 
examples of special Lagrangian tori in certain Calabi-Yau manifolds
that occur as hypersurfaces in complex projective space.
All of these are constructed as real slices of smooth
hypersurfaces defined over the reals.  This method of
constructing special Lagrangian submanifolds is, of course, well known.  
What does not appear to be in the current mirror symmetry literature 
is an explicit description of such examples in which the 
special Lagrangian submanifold is a 3-torus.

\medskip
{\it Some background.}
Let $M$ be a Ricci-flat \Kahler{} $n$-manifold, with \Kahler{} 
form~$\omega \in\cA^{1,1}_+(M)$ and a parallel holomorphic volume 
form~$\Omega\in\cA^{n,0}(M)$, normalized so that 
$$
\Omega\w\overline{\Omega}
= {{2^n\,(-i)^{n^2}}\over{n!}} \omega^n\,.
$$
The data~$(\omega,\Omega)$ will be said to constitute a {\it Calabi-Yau
structure} on~$M$.  
  
Given any oriented $\omega$-Lagrangian submanifold~$f:L\hookrightarrow M$
with induced volume form~$dV_L$, there exists a smooth 
map~$\lambda: L\to S^1\subset\bbC$ so that
$$
f^*(\Omega) = \lambda\,dV_L\,.
$$
The submanifold~$L$ is said to be {\it special Lagrangian\/}%
\footnote{$^1$}{Some authors require only that $\lambda$ be constant.  
For simplicity, I will not consider this (slight) generalization, 
though it is important for some purposes.}
if~$\lambda$ is identically~$1$.

As Harvey and Lawson~[HL] show, a compact special Lagrangian 
submanifold~$L\subset M$ has minimal $\omega$-volume in its 
homology class.  In fact, the real part of~$\Omega$ is a calibration 
on~$M$ and it calibrates the special Lagrangian submanifolds.
Consequently, any smooth (or even $C^1$) special Lagrangian 
submanifold is real analytic in~$M$.

McLean~[Mc] showed that for a smooth special
Lagrangian submanifold~$L\subset M$, the moduli space of nearby
special Lagrangian submanifolds is smooth and of dimension~$b_1(L)$.
He gave the following description:  Let~$\cH^1(L)\subset\cA^1(L)$
denote the vector space of harmonic 1-forms on~$L$ (with respect 
to the induced metric).  Then there exists an open neighborhood~$\cU$
of~$0\in\cH^1(L)$ and a
smooth (in fact, real analytic) mapping~$Q:\cU\to\cA^1(L)$
that vanishes to second order at~$0$ so that, in a~$C^1$-neighborhood 
of~$f:L\hookrightarrow M$, any special Lagrangian normal graph is 
of the form~$g_\alpha:L\to M$ where
$$
g_\alpha(x) 
=\exp_{f(x)}\bigl(i\,
f'\bigl(\bigl(\alpha{+}Q(\alpha)\bigr)^\sharp_x\bigr)\bigr),
$$
for a unique~$\alpha\in\cH^1(L)$.  (Here,  $\beta^\sharp$ denotes
the vector in~$TL$ dual to~$\beta\in T^*L$ and multiplication by~$i$
in~$TM$ carries vectors tangent to~$f(L)$ to vectors normal to~$f(L)$.)

In particular, if~$L$ is an $n$-torus and the induced metric on~$L$
is such that there is a basis~$\alpha_1,\ldots,\alpha_n$ of~$\cH^1(L)$ 
so that $\alpha_1\w\cdots\w\alpha_n$ is nowhere vanishing, 
then these normal graphs foliate a neighborhood of~$L$ in~$M$.
Such a special Lagrangian torus is said to be {\it foliating\/}. 
Every special Lagrangian $n$-torus is foliating when $n=1$~or~$2$,
but this is not known to be the case in higher dimensions.%
\footnote{$^{2}$}{This is trivial for~$n=1$.  For~$n=2$ this
follows from the fact a basis for the harmonic forms on~$L$
is formed by the real and imaginary parts of a nonzero (and
hence nonvanishing) holomorphic 1-form on~$L$ (where~$L$ is 
regarded as a complex curve of genus~$1$
via the conformal structure and orientation induced on~$L$ by its
inclusion into~$M$ as a special Lagrangian submanifold).} 

Conversely, if an open set~$U\subset M$ has a foliation~$\cL$ by compact
smooth special Lagrangian manifolds, then each leaf~$L$ of~$\cL$ is 
diffeomorphic to the $n$-torus and has a basis~$\alpha_1,\ldots,\alpha_n$ 
of~$\cH^1(L)$ so that $\alpha_1\w\cdots\w\alpha_n$ is nowhere vanishing.

In the last few years, these special Lagrangian foliations have become 
interesting in the theory of mirror symmetry, particularly after the 
paper of Strominger, Yau, and Zaslow~[SYZ].  Further work has been 
done on the geometry of moduli spaces of special Lagrangian 
tori and their augmentations by Hitchin~[Hi1,2] and Lu~[Lu], among others. 
 For example, in [Hi2], on sees how special Lagrangian foliations can
be constructed entirely in the holomorphic category in cases
where~$n$ is even and~$(M,\omega,\Omega)$ is actually hyper\Kahler.  

In general, explicit examples of special Lagrangian tori seem to be 
difficult to construct.  In the first place, the construction of Ricci-flat
\Kahler{} metrics itself relies on the celebrated theorem of Yau~[Ya],
which uses transcendental methods from nonlinear partial differential
equations to construct a unique Ricci-flat \Kahler{} form~$\omega$
in the cohomology class of any \Kahler{} form on a compact complex~$M$
that is endowed with a holomorphic volume form~$\Omega$.  (By 
a Bochner argument, $\Omega$ is then parallel with respect to~$\omega$.)

Many examples of compact complex \Kahler{} manifolds with a
holomorphic volume form are known, so Yau's theorem provides
a rich source of examples of Ricci-flat \Kahler{} manifolds.  
Perhaps the simplest examples (other than complex tori) are 
the smooth hypersurfaces of degree~$n{+}1$ in~$\bbC\bbP^n$.  

Given the data~$(M,\omega,\Omega)$ it is usually a highly nontrivial
matter to construct a special Lagrangian submanifold in~$M$.
One simple case where this can be done is
when~$(M,\omega,\Omega)$ possesses a {\it real structure\/}~$c$, 
i.e., an involution~$c:M\to M$, that satisfies~$c^*(\omega)=-\omega$ 
and~$c^*\Omega=\overline{\Omega}$.  If the locus~$L_c\subset M$
of fixed points of~$c$ is nonempty, then it is a smooth manifold of real
dimension~$n$ that can be oriented so as to be a special Lagrangian
submanifold of~$M$. Moreover, since~$c$ will necessarily preserve the 
underlying \Kahler{} metric, $L_c$ is totally geodesic in~$M$ with 
respect to this metric.

Of particular interest is the case of hypersurfaces~$H\subset\bbC\bbP^n$
of degree~$n{+}1$ that are invariant under the anti-automorphism induced by 
a real structure~$c:\bbC^{n+1}\to\bbC^{n+1}$, i.e., a complex antilinear
involution, with fixed subspace~$\bbR^{n+1}\subset\bbC^{n+1}$.  
In this case, the holomorphic volume form~$\Omega$ on~$H$ can be 
chosen to be anti\"\i{}nvariant under the induced 
antiholomorphism~$c_H:H\to H$, i.e., 
so that~$c_H^*(\Omega)=\overline{\Omega}$.
If one chooses a class~$\gamma\in H^{1,1}(H)$ that is representable
by a \Kahler{} form that also is anti\"\i{}nvariant under this conjugation,%
\footnote{$^3$}{These always exist. For example, take~$\gamma$ to be 
the class of the pullback to~$H$ of a $c$-invariant Fubini-Study metric 
on~$\bbC\bbP^n$.}
then the uniqueness part of Yau's theorem implies that the unique
Ricci-flat \Kahler{} form~$\omega$ on~$M$ representing the 
class~$\gamma$ must also be anti\"\i{}nvariant under~$c_H$.  
Consequently~$c_H$ is a real structure on~$(H,\omega,\Omega)$, 
implying that~$L_c=H\cap\bbR\bbP^n$ is a special Lagrangian submanifold.

Thus, one method of exhibiting a special Lagrangian torus is to find a
$c$-invariant smooth hypersurface~$H\subset\bbC\bbP^n$ whose real locus
$H\cap\bbR\bbP^n$ contains an $n$-torus as one of its components. 
In the following sections, I will point out some simple examples of this 
for~$n=2$,~$3$,~and~$4$.   Of course, the case~$n=2$ is well-known
and completely classical. Examples with~$n=3$ are almost as 
classical.  In fact, the possible topology of the real locus of a 
quartic surface in projective 3-space was determined by Kharlamov~[Kh], 
and there are several cases where this real locus contains a 2-torus
component.  Despite their simplicity, examples when~$n=4$ seem  
not to have been pointed out in the mirror symmetry literature before. 
Of course it is this case that is of the most interest. 

\medskip
{\it The construction.}
The basic idea is straightforward.  For any 
polynomial function~$P:\bbC^{n+1}\to \bbC$ 
that is homogeneous of degree~$d\ge1$, the equation $P=0$ defines
a cone~$C_P\subset\bbC^{n+1}$ and the equations~$P=dP=0$ define a 
subcone~$S_P\subset C_P$. Assume that $S_P$ is a proper subset of~$C_P$. 
Let~$\pi:\bbC^{n+1}\setminus\{{\bf0}\}\to\bbC\bbP^n$ be the usual
projection. Then the image~$H_P=\pi\bigl(C_P\setminus\{{\bf0}\}\bigr)$ 
is a hypersurface in~$\bbC\bbP^n$ that is smooth away from the 
singular locus~$\Sigma_P=\pi\bigl(S_P\setminus\{{\bf0}\}\bigr)$.

Let~$\cP$ denote the vector space consisting of polynomial 
functions~$P:\bbC^{n+1}\to\bbC$ that are homo\-geneous of degree~$n{+}1$ 
and that satisfy $P\circ c = \overline{P}$ where~$c:\bbC^{n+1}\to\bbC^{n+1}$ 
is the usual conjugation fixing~$\bbR^{n+1}\subset\bbC^{n+1}$.  This is
a real vector space of dimension~${{2n+1}\choose{n}}$.  For~$P\in\cP$,
the loci~$H_P$ and~$\Sigma_P$ are invariant under the induced conjugation 
on~$\bbC\bbP^n$.

Let~$\cO\subset\cP$ be the dense open subset consisting
of those~$P$ for which~$\Sigma_P$ is empty.  For~$P\in\cO$, endow~$H_P$ 
with the unique Ricci-flat \Kahler{} form~$\omega_P$ whose \Kahler{} 
class is dual to a hyperplane section.  Then~$\omega_P$ is 
anti\"{\i}nvariant under conjugation and the holomorphic volume
form~$\Omega_P$ can be chosen so that~$c^*(\Omega_P) = \bar\Omega_P$
and so that it satisfies the volume normalization needed to make
its real part have comass one.  This determines $\Omega_P$ up
to a sign.

Let~$\cT\subset\cP$ denote the open subset consisting of those~$P\in\cP$ 
for which the real slice~$H_P\cap\bbR\bbP^n$ contains 
an~$(n{-}1)$-torus~$L\subset H_P\cap\bbR\bbP^n$ that is disjoint 
from~$\Sigma_P$.  In the following sections, I will show that~$\cT$ 
is nonempty for $n=2$, $3$, and $4$ by producing explicit examples.
I do not know whether or not~$\cT$ is nonempty for higher~$n$.

When~$\cT$ is nonempty, it follows that $\cO\cap\cT$ is nonempty 
(and open).  For ~$P\in\cO\cap\cT$, the smooth 
hypersurface~$H_P\subset\bbC\bbP^n$ contains an $(n{-}1)$-torus as a 
component of its real locus.  Such an $(n{-}1)$-torus is special 
Lagrangian in the Calabi-Yau structure~$(H_P,\omega_P,\Omega_P)$.

\bigskip\bigskip
\centerline{\bf \S1. Cubic Curves in $\bbC\bbP^2$ }
\nobreak\bigskip\nobreak
This section is included for the sake of completeness and for
comparison with the cases where $n=3$ and $4$.
Any smooth cubic curve defined over~$\bbR$ is projectively equivalent 
to~$H_P$ where
$$
P = {X_0}^3 + {X_1}^3 + {X_2}^3 - 3\sigma\,X_0X_1X_2
$$
for some real number~$\sigma\not=1$.  The curve~$H_P$ has two real 
components when~$\sigma>1$ and one real component when~$\sigma<1$. 
In either case, there is always exactly one odd component%
\footnote{$^4$}{A simple closed curve in~$\bbR\bbP^2$ is {\it odd\/}
if it generates~$H_1(\bbR\bbP^2,\bbZ_2)\simeq\bbZ_2$, otherwise it
is {\it even}.}
and it contains the three real flexes~$f_1=[0,1,-1]$, 
$f_2=[-1,0,1]$, and $f_3=[1,-1,0]$.

When~$\sigma=1$, the curve~$H_P$ is singular, being the union 
of three lines, one of which is real, namely~$H_L$, where~$L=X_0+X_1+X_2$.
Note that the singular locus~$\Sigma_P$ consists of three points,
one of which is real, but which does not lie on~$H_L$.  

Thus, $P$ lies in~$\cT$ for any finite~$\sigma$.
(When~$\sigma=\infty$, the real locus consists of three real
non-concurrent lines and hence has no smooth component.)  Note also
that $\cO\cap\cT$ has two components.

\bigskip\bigskip
\centerline{\bf \S2. Quartic Surfaces in $\bbC\bbP^3$ }
\nobreak\bigskip\nobreak
The topology of the real locus of a quartic surface in projective 3-space 
was determined in the 1970s.  For a survey of these results, see
Kharlamov~[Kh].   However, it is easy to construct elements of~$\cT$
directly.  Here are a few simple examples, generalizing both the 
`odd component' and the `even component' case of curves.
\smallskip

As a first example, consider
$$
P = {X_0}^4 + {X_1}^4 - {X_2}^4 - {X_3}^4.
$$
This is a nonsingular Fermat-type quartic.  The real locus is
$$
H_P\cap\bbR\bbP^3 = 
\{ [a,\,b,\,u,\,v]\in\bbR\bbP^3\ \vrule\  a^4+b^4=u^4+v^4=1 \}.
$$
Since the curve~$C\subset\bbR^2$ defined by~$x^4+y^4=1$ is diffeomorphic
to the circle,~$H_P\cap\bbR\bbP^3$ is diffeomorphic
to~$(C\times C)/{\sim}$ 
where 
$\bigl((a,\,b),\,(u,\,v)\bigr)\sim\bigl(({-}a,\,{-}b),\,({-}u,\,{-}v)\bigr)$. 
This is manifestly a 2-torus.  Consequently, $P$ lies in~$\cT$.  

Note that $\bbR\bbP^3\cap H_P$ is like the odd component of an elliptic
curve: It is not contractible in~$\bbR\bbP^3$ because 
its inverse image under the double cover~$S^3\to\bbR\bbP^3$ 
is a nontrivial double cover.

\bigskip
As a second example, take
$$
R = \bigl({X_0}^2 + {X_1}^2 - r_1({X_2}^2 + {X_3}^2)\bigr)
    \bigl({X_0}^2 + {X_1}^2 - r_2({X_2}^2 + {X_3}^2)\bigr)
$$
where~$0<r_1<r_2$.  In this case, the real locus $\bbR\bbP^3\cap H_R$ 
consists of {\it two\/} disjoint 2-tori, neither of which is
contractible in~$\bbR\bbP^3$.  The singular locus~$\Sigma_R$ 
consists of the four lines
$$
L_{m,n} = \left\{\ [a,\,(-1)^m\,ia,\,b,\,(-1)^n\,ib]\>\vrule\> 
           [a,b]\in\bbC\bbP^1\ \right\},\qquad m,n\in\{0,1\}.
$$ 
and these have no real points.

\bigskip
The third example is more like the even component
of a real elliptic curve.  In~$\bbR^3$, consider the circle~$C$ 
defined by the equations
$$
{x_1}^2 + {x_2}^2 - 1 = x_3 = 0.
$$
This is the minimum locus of the quartic polynomial~$q_0$ defined by
$$
q_0 = \bigl({x_1}^2 + {x_2}^2 - 1\bigr)^2 + {x_3}^4.
$$
The critical locus of~$q_0$ in~$\bbR^3$ is the circle~$C$ plus
the origin~$O = (0,0,0)$.  Consequently, the regular values of~$q_0$
are all real numbers other than~$0$ and~$1$.  

Let~$\epsilon$ be a real number satisfying~$0<\epsilon<1$ and consider
the quartic polynomial~$q_\epsilon = q_0-\epsilon$.  Its zero locus
in~$\bbR^3$ is smooth and is the boundary of the 
region~$R_\epsilon = {q_0}^{-1}\bigl([0,\epsilon]\bigr)$, which 
retracts onto the circle~$C$.  Consequently, the zero locus of~$q_\epsilon$
is diffeomorphic to a torus.

Now consider the homogeneous quartic~$Q$ defined by
$$
\eqalign{
Q 
&= \bigl({X_1}^2+{X_2}^2-{X_0}^2\bigr)^2+{X_3}^4-\epsilon\,{X_0}^4\cr
&= {X_0}^4\,q_\epsilon\bigl(X_1/X_0,X_2/X_0,X_3/X_0\bigr)\,.
}
$$
The singular locus~$\Sigma_Q$ consists of two
nonreal points~$[0,1,\pm i,0]$.

Since~$Q=X_0=0$ on~$\bbR^4$ only at the origin, the real 
slice~$H_Q\cap\bbR\bbP^3$ is just
$$
T = \left\{\,[1,X_1,X_2,X_3]\ 
        \vrule\ q_\epsilon\bigl(X_1,X_2,X_3\bigr)=0 \,\right\}
$$
and hence is diffeomorphic to a torus.  Thus~$Q$ lies in~$\cT$. 

Note that~$T$ misses a linear~$\bbR\bbP^2$ in~$\bbR\bbP^3$ and
consequently is contractible in~$\bbR\bbP^3$.
Since~$H_Q\cap\bbR\bbP^3$ is not homotopic to~$H_P\cap\bbR\bbP^3$,
it follows that $P$ and $Q$ lie in different components of~$\cT$.
\bigskip
By these three examples, it follows that~$\cO\cap\cT$ has at least 
three components.

\bigskip\bigskip
\centerline{\bf \S3. Quintic $3$-folds in $\bbC\bbP^4$ }
\nobreak\bigskip\nobreak
In~$\bbR^4$, consider the locus~$T$ defined by the equations
$$
\eqalign{
{x_1}^2+{x_2}^2&=1, \cr
{x_3}^2+{x_4}^2&=1. \cr
}
$$
This is a 2-torus (it is the Clifford torus, up to scale). 
Because it is cut out by two independent equations, the normal bundle~$N$ 
of~$T$ in~$\bbR^4$ is trivial, implying that its unit circle 
bundle~$L\subset N$ is diffeomorphic to the $3$-torus.

Consider the quartic function
$$
q = \bigl({x_1}^2+{x_2}^2-1\bigr)^2
      + \bigl({x_3}^2+{x_4}^2-1\bigr)^2.
$$
Then $q(x)\ge0$, with equality if and only if~$x$ lies in~$T$.  
The critical values of~$q$ are $0$, $1$, and $2$. 

Thus, for~$\epsilon$ satisfying~$0<\epsilon<1$, the
hypersurface~$T_\epsilon\subset\bbR^4$ defined by~$q(x)-\epsilon^2=0$
is the boundary of a compact domain that has~$T$ as a deformation
retract.  In fact, by Morse theory, ~$T_\epsilon$ is diffeomorphic 
to~$L\subset N$, which, as has been remarked, is a $3$-torus.

Now consider the (reducible) quintic polynomial~$P$ defined by
$$
\eqalign{
P &= {X_0}^5\,\bigl(q(X_1/X_0,X_2/X_0,X_3/X_0,X_4/X_0) - \epsilon^2\bigr)\cr
&=X_0\,\bigl[\bigl({X_1}^2+{X_2}^2-{X_0}^2\bigr)^2
      + \bigl({X_3}^2+{X_4}^2-{X_0}^2\bigr)^2-\epsilon^2\,{X_0}^4\bigr]\cr
&= X_0\,Q,\qquad\hbox{($Q$ is irreducible).}\cr
}
$$

The quintic hypersurface~$H_P\subset\bbC\bbP^4$ is the union of the 
hyperplane~$H_{X_0}$ and the quartic hypersurface~$H_Q$.

The singular locus~$\Sigma_Q$ is a union of four lines
$$
L_{m,n} = \left\{\ [0,\,a,\,(-1)^m\,ia,\,b,\,(-1)^n\,ib]\>\vrule\> 
           [a,b]\in\bbC\bbP^1\ \right\},\qquad m,n\in\{0,1\}.
$$ 
None of these lines have real points.

The intersection~$H_{X_0}\cap H_Q$ is a union of two quadric surfaces
$$
S_{\pm} = \left\{\ [0,X_1,X_2,X_3,X_4]\>\vrule\> 
           \bigl({X_1}^2+{X_2}^2\bigr) \pm i\bigl({X_3}^2+{X_4}^2\bigr) = 0
            \ \right\}.
$$
Neither of these surfaces has any real points.  Note that $S_+\cap S_-$
consists of the four lines~$L_{m,n}$. Consequently,~$\Sigma_P=S_+\cup S_-$.

Thus, the real slice~$H_P\cap\bbR\bbP^4$ is smooth and
is the disjoint union of~$H_{X_0}\cap\bbR\bbP^4\simeq\bbR\bbP^3$
and~$H_Q\cap\bbR\bbP^4\simeq T_\epsilon$, which is diffeomorphic
to the $3$-torus.

Thus,~$P$ lies in~$\cT$, which is thereby shown to be nonempty.
Consequently,~$\cO\cap\cT$ is nonempty, thus proving that there
are smooth (and hence irreducible) quintic hypersurfaces in~$\bbC\bbP^4$ 
whose real locus contains a torus as a component.

\bigskip\bigskip
\centerline{\bf \S4. References }
\nobreak\bigskip\nobreak

\item{[Hi1]} 
{\smc N. Hitchin}, 
{\it The moduli space of special Lagrangian submanifolds},
dg-ga/9711002.
\smallskip

\item{[Hi2]} 
{\smc N. Hitchin}, 
{\it The moduli space of complex Lagrangian submanifolds},
dg-ga/9901069.
\smallskip

\item{[HL]} 
{\smc F. Harvey} and {\smc H. Lawson}, 
{\it Calibrated Geometries}, Acta. Math.~{\bf 148} (1982), 47--157. 
MR~85i:53058.
\smallskip

\item{[Kh]}
{\smc V. Kharlamov},
{\it On the classification of nonsingular surfaces of degree~$4$ 
in $\bbR\bbP^3$ with respect to rigid isotopies.} (Russian) 
Funktsional. Anal. i Prilozhen. {\bf 18} (1984), 49--56.
(English translation in Functional Anal. Appl. {\bf 18} (1984), 39--45.)
MR 85m:14034.
\smallskip

\item{[Lu]}
{\smc P. Lu}, {\it Special Lagrangian Tori on a Borcea-Voisin Threefold\/},
dg-ga/9902063.
\smallskip

\item{[Mc]} 
{\smc R. McLean}, 
{\it Deformations of calibrated submanifolds}, 
Comm. Anal. Geom.Ê {\bf 6} (1998), 705--747. 
\smallskip

\item{[SYZ]} 
{\smc A. Strominger, S.T. Yau,} and {\smc E. Zaslow},
{\it Mirror Symmetry is T-Duality}, Nucl. Phys. {\bf B479} (1996),
243--259. MR 97j:32022.
\smallskip

\item{[Ya]} 
{\smc S.T. Yau}, 
{\it On the Ricci curvature of a compact K\"{a}hler manifold 
and the complex Monge-Amp\'{e}re equations. I}, 
Comm. Pure Appl. Math. ~{\bf 31} (1978), 339--411.  MR 81d:53045.

\bye